\documentclass{article}

\usepackage{fullpage}   
\usepackage{natbib}

\usepackage[utf8]{inputenc} % allow utf-8 input
\usepackage[T1]{fontenc}    % use 8-bit T1 fonts
\usepackage{hyperref}       % hyperlinks
\usepackage{url}            % simple URL typesetting
\usepackage{booktabs}       % professional-quality tables
\usepackage{amsfonts}       % blackboard math symbols
\usepackage{nicefrac}       % compact symbols for 1/2, etc.
\usepackage{microtype}      % microtypography
\usepackage{xcolor}         % colors

\usepackage{nicematrix}

\usepackage{algorithm}
\usepackage{algorithmic}

\newcommand{\cD}{\mathcal{D}}

\newcommand{\cF}{\mathcal{F}}

\newcommand{\cN}{\mathcal{N}}

\newcommand{\cU}{\mathcal{U}}

\newcommand{\cX}{\mathcal{X}}

\newcommand{\cW}{\mathcal{W}}
\usepackage{mathtools}
\newcommand{\eqdef}{\coloneqq}
\newcommand{\R}{\mathbb{R}} 

\newcommand{\Exp}[1]{{\rm E} \left[ #1 \right]} 
\newcommand{\ExpCond}[2]{{\rm E}\left[#1\;|\;#2\right]}
\newcommand{\ExpSub}[2]{{\rm E}_{#1}\left[#2\right]}
\newcommand{\sqn}[1]{{\left\| #1 \right\|^2}} % squared norm (alternative)
\newcommand{\abr}[1]{\left\langle#1\right\rangle} % angle bracket

\usepackage{color}
\definecolor{ForestGreen}{RGB}{34,139,34}
\usepackage{xspace}
\usepackage{scalefnt}
\newcommand{\algname}[1]{{\sf\color{ForestGreen}\scalefont{0.95}{#1}}\xspace}

\usepackage{amsmath, amsfonts,amssymb,amsthm}
\newtheorem{theorem}{Theorem}

\newtheorem{lemma}{Lemma}
\newtheorem{assumption}{Assumption}

\title{\textbf{Stochastic Proximal Point Methods\\ for Monotone Inclusions under Expected Similarity}}

 \author{Abdurakhmon Sadiev$^{1}$ \qquad  Laurent Condat$^{1,2}$  \qquad Peter Richt\'{a}rik$^{1,2}$\\
 \phantom{xx}
 \\
$^1$Computer Science Program, CEMSE Division,\\ King Abdullah University of Science and Technology (KAUST)\\ Thuwal, 23955-6900, Kingdom of Saudi Arabia\\
 $^2$SDAIA-KAUST Center of Excellence in Data Science and \\Artificial Intelligence 
(SDAIA-KAUST AI)
}

\date{May 22, 2024}

\begin{document}

\maketitle

\begin{abstract}
Monotone inclusions have a wide range of applications, including minimization, saddle-point, and equilibria problems. We introduce new stochastic algorithms, with or without variance reduction, to estimate a root of the expectation of possibly set-valued monotone operators, using at every iteration one call to the resolvent of a randomly sampled operator. We also introduce a notion of similarity between the operators, which holds even for discontinuous operators.
We leverage it to derive linear convergence results in the strongly monotone setting.
\end{abstract}

\tableofcontents

\section{Introduction}

We consider stochastic monotone inclusions in a given finite-dimensional real Hilbert space $\cX$, 
 which are problems of the form 
\begin{equation}
    \label{eq:stoch_mon_inc_problem}
    \text{Find}~~x^{\star} \in \cX~~\text{such that}~~0\in A(x^{\star}),~~\text{where}~~A\eqdef \ExpSub{\xi\sim\cD}{A_{\xi}}
\end{equation}
and $A_{\xi}$ is a possibly set-valued monotone operator for every random sample $\xi$ of a distribution $\cD$. We recall basic  notions of monotone operator theory in Section \ref{secmo} and refer to the  textbook \citet{bau17} for more details. For instance, when $\cD$ is the uniform distribution over $[n]\eqdef\{1,\ldots,n\}$ for some $n\geq 2$,  \eqref{eq:stoch_mon_inc_problem} becomes the finite-sum monotone inclusion
\begin{equation}
    \label{eq:finite_sum_mon_inc_problem}
    \text{Find}~~x^{\star} \in \cX~~\text{such that}~~0\in A(x^{\star})\eqdef \frac{1}{n}\sum^n_{i=1}A_i(x^{\star}).
\end{equation}
We introduce randomized algorithms, with or without variance reduction, to solve \eqref{eq:stoch_mon_inc_problem}. They use at every iteration the resolvent of one randomly chosen $A_{\xi}$.

\subsection{Motivation}

Monotone inclusions \citep{sim08,bau17} have a wide range of applications \citep{kap01,fac03,com13}, 
in 
mechanics \citep{glo84,glo89}, partial differential equations \citep{mer793,att08,gho09,par24}, mean field games \citep{bri18,gom21}, 
control \citep{sta16}, 
communications \citep{pal09,xin20}, 
traffic equilibrium \citep{fuk96,att10},
optimal transport \citep{pap14}, 
Nash equilibria and game theory \citep{mor53, bri13,bra18,yi19,loi21}, 
and are of utmost importance in machine learning. 
 Primarily, they encompass optimization problems \citep{sra11,bac12,cev14,pol15, bub15}:
minimizing a convex function $f: \cX\rightarrow \R\cup\{+\infty\}$ is equivalent to \eqref{eq:stoch_mon_inc_problem} with $A=\partial f$, the subdifferential of $f$, and finding a stationary point of a smooth but possibly nonconvex function $f$ is equivalent to \eqref{eq:stoch_mon_inc_problem} with $A=\nabla f$, the gradient of $f$. Nonconvex nonsmooth optimization problems have variational formulations, too \citep{mor06}. 
Moreover, splitting algorithms to solve structured optimization problems can be derived by formulating the problem as a monotone inclusion in a higher-dimensional lifted space. For instance, minimizing $f+\sum_{i=1}^n g_i(K_i x)$, for linear operators $\cW\rightarrow \cU_i$ and functions $g$ and $h_i$, can be formulated as \eqref{eq:stoch_mon_inc_problem} with $\cX=\cW\times \cU_1\times\cdots\times\cU_n$ and the monotone operator 
\begin{equation}
A = \left(\begin{array}{cccc}\partial{f}& K_1^* u_1&\cdots& K_n^* u_n\\
-K_1x&(\partial g_1)^{-1}(u_1)&0&0\\
\vdots&0&\ddots&0\\
-K_nx&0&0&(\partial g_n)^{-1}(u_n)
\end{array}\right),
\end{equation}
where $\cdot^*$ denotes the adjoint operator. For a suitable preconditioning linear operator $P$, that is symmetric and positive definite, $P^{-1}A$ is monotone in $\cX$ endowed with the modified inner product $\langle \cdot, P\cdot\rangle$, and one can design iterative algorithms to solve $0\in P^{-1}A(x)$ \citep{bau08, he10, bri11, com12,con13,com14, com18,sal20,com21,con22,sal22,con23}.

Besides minimization problems, monotone inclusions allow us to formulate saddle-point problems \citep{leb67,roc64,mcl74,nem04,dav162,bui22,dav23}, which have many applications in machine learning  \citep{mer19,gor22,bez23}, e.g.\ for adversarial training \citep{goo15,mad18}, GANs \citep{gid19}, and distributionally robust optimization \citep{nam16}.

We propose different algorithms in the framework of the Stochastic Proximal Point Method (SPPM), with or without variance reduction. Even in the optimization setting, our study under a similarity assumption, which is weaker than smoothness, is new, to the best of our knowledge.

\section{Definitions and Properties of Monotone Operators}\label{secmo}

Let $B:\cX\rightarrow 2^{\cX}$ be a set-valued operator on $\cX$. We define its graph $\mathrm{gra}(B)\eqdef\{ (x,u)\in\cX^2\;:\;u\in B(x)\}$ and its inverse $B^{-1}$ as the set-valued operator whose graph is  $\mathrm{gra}(B^{-1})\eqdef\{ (u,x)\in\cX^2\;:\;u\in B(x)\}$. $x\in\cX$ is a \emph{zero} of $B$ if  $0\in B(x)$. 

\subsection{Monotone Operators}

$B$ is \emph{monotone} if for every $(x,u)$ and $(y,v)$ in $\mathrm{gra}(B)$, 
\begin{equation*}
\langle u-v,x-y\rangle \geq 0.
\end{equation*}
$B$ is \emph{maximally monotone} if there exists no monotone operator whose graph strictly contains 
$\mathrm{gra}(B)$. $B$ is (maximally) monotone if and only if $B^{-1}$ is (maximally) monotone. The subdifferential $\partial f$ of a proper lower semicontinuous convex function $f$ is maximally monotone. 

$B$ is \emph{$\mu$-strongly monotone} for some $\mu>0$ if, for every $(x,u)$ and $(y,v)$ in $\mathrm{gra}(B)$, 
    \begin{equation}
        \label{eq:strong_mon}
        \abr{u-v, x-y } \geq \mu\sqn{x-y}.
    \end{equation}
In that case, $\gamma B$ is $\gamma\mu$-strongly monotone, for every $\gamma>0$. If $B$ is $\mu$-strongly maximally monotone, its zero exists and is unique.

The following assumption on the operators in    \eqref{eq:stoch_mon_inc_problem} will be considered to analyze the proposed algorithms. 
\begin{assumption}[strong monotonicity]
\label{as:strong_mon}
    There exists $\mu>0$ such that $A_{\xi}$  is $\mu$-strongly maximally monotone for every $\xi\sim\cD$. 
        Therefore, $A\eqdef \ExpSub{\xi\sim\cD}{A_{\xi}}$  is $\mu$-strongly maximally monotone as well and 
    the solution $x^\star$ to \eqref{eq:stoch_mon_inc_problem} exists and is unique.
\end{assumption}

A single valued operator $C:\cX\rightarrow \cX$ is $\beta$-\emph{cocoercive} for some $\beta>0$ if, for every $(x,y)\in\cX^2$, 
    \begin{equation*}
       \left\langle x-y,C(x)-C(y) \right\rangle \geq \beta\sqn{C(x)-C(y)}.
    \end{equation*}
%bau17, Example 20.31
A function $f$ is \emph{$L$-smooth} for some $L>0$ if it is differentiable and its gradient $\nabla f$ is $L$-Lipschitz continuous;  that is, for every $(x,y)\in\cX^2$, 
    \begin{equation*}
     \|\nabla f(x)-\nabla f(y)\|\leq L \|x-y\|.
     \end{equation*}
In that case, $\nabla f$ is $L^{-1}$-cocoercive, according to the Baillon--Haddad theorem. This equivalence between Lipschitz-continuity and cocoercivity only holds for operators which are gradients of convex functions. In general, a monotone operator can be Lipschitz-continuous without being cocoercive. A prominent example is the skew operator $(x,y)\in\cX^2\mapsto (K^*y, Kx)$ for any linear operator $K$ on $\cX$, which is $\|K\|$-Lipschitz continuous but not cocoercive. Thus, monotone inclusions are much more general than optimization problems. In particular, the forward algorithm generalizing gradient descent, which iterates $x^{k+1} \eqdef x^k - \gamma A(x^k)$ for a maximally monotone single-valued operator $A$ and a stepsize $\gamma>0$, converges if $A$ is cocoercive, but not if it is merely Lipschitz-continuous (take $-I$ as an example, where  $I$ denotes the identity: the iteration diverges for every $\gamma>0$ and $x^0\neq 0$). This is why robust iterative fixed-point algorithms to solve monotone inclusions use the resolvent of the monotone operators, as we describe in the next section.

\subsection{The Resolvent and the Proximal Point Method}

The \emph{resolvent} of $B$ is the operator $(I + B )^{-1}$. 
According to the Minty theorem, if $B$ is maximally monotone, its resolvent is defined everywhere and single-valued. 
The resolvent of a strongly monotone operator is contractive:
\begin{lemma}[contractivity of the resolvent]
    \label{fact:non-expansiveness}
    Let $B:\cX\rightarrow 2^{\cX}$ be a $\mu$-strongly maximally monotone operator, for some $\mu>0$. Then its resolvent is $(1+\mu)^{-1}$-contractive; that is, for every $(x,y)\in\cX$, 
    \begin{equation}
        \label{eq:non-expansiveness}
        \|x^+-y^+\| \leq \frac{1}{1+\mu}\|x-y\|,
    \end{equation}
    where $x^+ = (I+B)^{-1}(x)$, $y^+ = (I+B)^{-1}(y)$.
\end{lemma}
The resolvent of the subdifferential $\partial f$ of a function $f$ is its proximity operator  $\mathrm{prox}_f = (I+\partial f)^{-1}: x\in\cX \mapsto \arg\min_y \big(f(y) + \frac{1}{2}\|y-x\|^2\big)$. Optimization algorithms making use of proximity operators are called proximal algorithms \citep{par14,ryu16,con23}. The iteration $x^{k+1}\eqdef \mathrm{prox}_f (x^k)$ to minimize a function $f$, and by extension the iteration $x^{k+1}\eqdef (I+B)^{-1}(x^k)$ to find a zero of the operator $B$, is called the \emph{proximal point algorithm}, or \emph{proximal point method} (PPM) \citep{roc76}. It follows from Lemma~ \ref{fact:non-expansiveness} that if $B$ is $\mu$-strongly maximally monotone, the PPM converges linearly to its zero $x^\star=(I+B)^{-1}(x^\star)$, which exists and is unique, since $\|x^{k+1}-x^\star\| \leq \frac{1}{1+\mu}\|x^k-x^\star\|$ for every $k\geq 0$.

\subsection{Similarity Between Operators}

It is natural to consider that there  exists some level of similarity or homogeneity between the operators $A_i$, in particular in machine learning where they express characteristics of underlying data \citep{hen20,sun22,cha24}. 
To capture this property, we define two notions of similarity.
\begin{assumption}[expected similarity]
    \label{as:similarity_operator}
    In Problem~\eqref{eq:stoch_mon_inc_problem}, there exist a solution $x^\star$ and $\delta >0$ such that, for every $x \in \cX$ and $a_\xi \in A_{\xi}(x)$, $\xi\sim\cD$, 
    there exist $a_{\xi}^{\star}\in A_{\xi}(x^{\star})$, $\xi\sim\cD$, such that $\ExpSub{\xi\sim\cD}{a_{\xi}^{\star}}=0$, and 
    \begin{equation}
        \label{eq:similarity_operator}
          \ExpSub{\xi\sim\cD}{\sqn{a_{\xi} - \ExpSub{\xi'\sim\cD}{a_{\xi'}} - a_{\xi}^{\star}}} \leq \delta^2\sqn{x-x^{\star}}.
    \end{equation}
\end{assumption}
This assumption can be satisfied by set-valued operators with discontinuities and is even weaker than assuming every $A_\xi-A$ to be Lipschitz-continuous at $x^\star$.  
 Let us give a simple example of $n=2$ maximally monotone operators $A_1: x\in\mathbb{R} \mapsto (\{1\}$ if $x<1$, $[1,3]$ if $x=1$, $\{3\}$ if $x>1)$ and $A_2: x\in\mathbb{R} \mapsto (\{4x-7\}$ if $x<1$, $[-3,-1]$ if $x=1$, $\{4x-5\}$ if $x>1)$ on $\cX=\mathbb{R}$, with $\cD$ the uniform distribution on $[n]$. We have $A=\frac{1}{2}(A_1+A_2):x\in\mathbb{R} \mapsto (\{2x-3\}$ if $x<1$, $[-1,1]$ if $x=1$, $\{2x-1\}$ if $x>1)$, and $x^\star=1$. For every $x<1$, with $a_1=1\in A_1(x)$, $a_2=4x-7\in A_2(x)$, $a_1^\star=2$,  $a_2^\star=-2$, we can check that 
 \eqref{eq:similarity_operator} is satisfied with $\delta=2$, as the left-hand side is $(2x-2)^2$. For every $x>1$, with $a_1=3\in A_1(x)$, $a_2=4x-5\in A_2(x)$, $a_1^\star=2$,  $a_2^\star=-2$, we can check that 
 \eqref{eq:similarity_operator} is satisfied with $\delta=2$, as the left-hand side is $(2x-2)^2$ as well. At $x=x^\star=1$, for every $a_1\in [1,3]= A_1(x)$ and $a_2\in [-3,-1]= A_2(x)$, with $a_1^\star = \frac{1}{2}(a_1-a_2) =-a_2^\star$,  \eqref{eq:similarity_operator} is satisfied with any $\delta$, as the left-hand side is zero. Overall, \eqref{eq:similarity_operator} is satisfied $\delta=2$.\smallskip

\begin{assumption}[average similarity]
    \label{as:another_similarity_operator}
    In Problem~\eqref{eq:finite_sum_mon_inc_problem}, there exist a solution $x^\star$ and $\tilde{\delta} >0$ such that, for every $x_i \in \cX$  and $a_i \in A_i(x_i)$, $i\in[n]$,     there exist $a_i^{\star}\in A_i(x^{\star})$, $i\in[n]$, such that $\sum_{i=1}^n a_i^\star = 0$, and 
    \begin{equation}
        \label{eq:another_similarity_operator}
        \frac{1}{n}\sum_{i=1}^n\sqn{a_i - \frac{1}{n}\sum_{j=1}^n a_j  - a_i^{\star}} \leq  \frac{\tilde{\delta}^2}{n}\sum_{i=1}^n\sqn{x_i-x^{\star}}.
    \end{equation}
\end{assumption}
Assumption~\ref{as:another_similarity_operator} is stronger than Assumption~\ref{as:similarity_operator}  with $\cD$ the uniform distribution over $[n]$, since \eqref{eq:another_similarity_operator} with $x_1=\cdots=x_n=x$ implies \eqref{eq:similarity_operator}.

Related definitions of similarities have been considered in several works \citep{hen20,sun22,szl22,kov22,kha23,lin23}. 
For instance, the property that
for every $(x,y) \in \cX^2$
    \begin{equation*}
        %\label{eq:av_similarity}
        \frac{1}{n}\sum_{i=1}^n\sqn{A_i(x) - A(x) - A_i(y) + A(y)} \leq \delta^2  \sqn{x-y},
    \end{equation*}
in the case where the $A_i=\nabla f_i$ are gradients of smooth functions $f_i$, is called 
Hessian variance in \citet{szl22} and $\delta$-average second-order similarity in \citet{lin23}.  
Indeed, if the functions $f_i$ are twice differentiable, this property is equivalent to the one that, for every $x \in \cX$,     
\begin{equation*}
        %\label{eq:av_similarityhess}
        \frac{1}{n}\sum_{i=1}^n\sqn{\nabla^2 f_i(x)-\nabla^2 f(x)} \leq \delta^2;
    \end{equation*}
    that is, the variance of the Hessians $\nabla^2 f_i$ is uniformly bounded.

\section{The Stochastic Proximal Point Method (SPPM)}

\begin{algorithm}[H]
    \begin{algorithmic}[1]
        \STATE  \textbf{Parameters:}  stepsize $\gamma>0$, initial estimate $x^0\in\cX$
        \FOR {$k=0,1, \ldots$}
        \STATE Sample $\xi^k\sim \cD$
        \STATE $x^{k+1} \eqdef \left(I+\gamma A_{\xi^k}\right)^{-1}(x^k)$
        \ENDFOR
    \end{algorithmic}
    \caption{Stochastic Proximal Point Method (\algname{SPPM})}
    \label{alg:sppm}    
\end{algorithm}

The Stochastic Proximal Point Method (SPPM) consists of iterating the resolvent of an operator $A_{\xi^k}$ chosen randomly at every iteration $k$. Under Assumption~\ref{as:strong_mon}, it converges linearly to a neighborhood of $x^\star$.

\begin{theorem}
    \label{th:sppm_convergence}
   In  Problem \eqref{eq:stoch_mon_inc_problem}, let Assumption~\ref{as:strong_mon} hold, 
   and for every $\xi\sim\cD$, let $a_\xi^\star \in A_{\xi}(x^{\star})$, such  that $\ExpSub{\xi\sim\cD}{a_\xi^\star}=0$. Such $a_\xi^\star$ exist by definition of $x^{\star}$. If they are not unique, we define them as ones minimizing 
            \begin{equation}
        \label{eq:sigma_star}
        \sigma^2_{\star} \eqdef \ExpSub{\xi\sim\cD}{\sqn{a_\xi^\star}}.
    \end{equation}
    %$\ExpSub{\xi\sim\cD}{\sqn{a_\xi^\star}}$.
    Then in \algname{SPPM} with any stepsize $\gamma>0$ and initial estimate $x^0\in\cX$, we have,   for every $k \geq 0$,
    \begin{align}
        \Exp{\sqn{x^{k}-x^{\star}}} &\leq  \left(\frac{1}{1+\gamma\mu}\right)^{2k} \sqn{x^0 - x^{\star}} + \frac{1-(1+\gamma\mu)^{-2k}}{(1+\gamma\mu)^2-1}\gamma^2\sigma^2_{\star}\label{eq9}\\
             &\leq \left(\frac{1}{1+\gamma\mu}\right)^{2k} \sqn{x^0 - x^{\star}} + \frac{\gamma\sigma^2_{\star}}{2\mu+\gamma\mu^2}.\label{eq10}
    \end{align}
    %where
\end{theorem}
Our result is tight: \eqref{eq9} is satisfied with an equality with the operators $A_{\xi}(x)=\mu(x-x^\star)+a_\xi^\star$ for some $\mu>0$, $x^\star\in \cX$, and  $a_\xi^\star \in \cX$ such  that $\ExpSub{\xi\sim\cD}{a_\xi^\star}=0$.

Even in the optimization setting,  
Theorem~\ref{th:sppm_convergence} is new. 
 In \citet{ber11}, the SPPM, called \emph{incremental proximal algorithm}, was studied  to minimize a finite sum of functions, but the convergence bounds depend on the number of functions, so they are not applicable to our setting where the distribution $\cD$ is arbitrary. In \citet{bia16} and \citet{tou16}, convergence results with decreasing stepsizes are derived. 
SPPM-type algorithms have been studied for stochastic optimization in \citet{asi19}, with a focus on stability in the case of inexact computation of the proximity operator. In \citet{dav19} the SPPM is studied for optimization, but their convergence analysis (Theorem 4.4) relies on the decay of the function values, so it is not applicable to our setting.

In \citet[Theorem 7]{ryu14}, in the convex optimization setting, by simply using the triangular inequality $\|x^{k+1}-x^\star\|\leq \|\left(I+\gamma A_{\xi^k}\right)^{-1}(x^k)-\left(I+\gamma A_{\xi^k}\right)^{-1}(x^\star)\|+\|\left(I+\gamma A_{\xi^k}\right)^{-1}(x^\star)-x^\star\|\leq (1+\gamma\mu)^{-1}\|x^k-x^\star\|+\gamma\|\tilde{a}_{\xi^k}^\star\|$, where $\tilde{a}_\xi^\star$ is the minimum-norm element of $A_{\xi}(x^\star)$, 
they obtain
    \begin{equation*}
        \Exp{\|x^{k}-x^{\star}\|} \leq \left(\frac{1}{1+\gamma\mu}\right)^{k}\|x^0-x^{\star}\| + \frac{(1+\gamma\mu)\tilde{\sigma}_\star}{\mu},
    \end{equation*}
    where   $\tilde{\sigma}_{\star} = \ExpSub{\xi\sim\cD}{\|\tilde{a}_\xi^\star\|}$.   The neighborhood size does not tend to zero when $\gamma\rightarrow 0$, as is the case in \eqref{eq10}. 
In \citet[Theorem 10]{pat18}, the following result is obtained in the convex optimization setting with \emph{smooth} functions:
    \begin{equation*}
        \Exp{\sqn{x^{k}-x^{\star}}} \leq 2\left(\frac{1}{1+\gamma\mu}\right)^{2k}\sqn{x^0-x^{\star}} + \frac{2(1+\gamma\mu)^2\sigma^2_{\star}}{\gamma^2}.
    \end{equation*}
    The neighborhood size tends to $+\infty$ when $\gamma\rightarrow 0$, whereas it should tend to zero. In the same setting, 
     \citet[eq.~19]{kha23} derived
    \begin{equation*}
        \Exp{\sqn{x^{k}-x^{\star}}} \leq \left(\frac{1}{1+\gamma\mu}\right)^{k}\sqn{x^0-x^{\star}} + \frac{\gamma\sigma^2_{\star}}{\mu}.
    \end{equation*}
The rate and the neighborhood size are larger than in \eqref{eq10}. Thus, even in the optimization setting, our result is new and tight, with a simple and elegant proof.

\section{SPPM with Operator Correction (SPPM-OC)}

The \algname{SPPM} does not converge to the exact solution $x^\star$ of \eqref{eq:stoch_mon_inc_problem} but only to its neighborhood. To correct this shortcoming, 
we propose a new algorithm, the SPPM with Operator Correction (\algname{SPPM-OC}). It is variance-reduced \citep{gow20a}; that is, it converges to the exact solution under Assumptions \ref{as:strong_mon} and \ref{as:similarity_operator}. This is achieved by adding a shift to $x^k$ before applying the resolvent of a randomly chosen $A_{\xi^k}$, to correct for the difference between $A_{\xi^k}$ and its expectation $A$.

\begin{algorithm}[H]
	\begin{algorithmic}[1]
		\STATE  \textbf{Parameters:} stepsize $\gamma>0$, initial estimate $x^0\in\cX$
		\FOR {$k=0,1, \ldots$}
		\STATE Sample $\xi^k\sim \cD$
		\STATE Choose $a^k_{\xi^k}\in A_{\xi^k}(x^k)$ and $a^k \in A(x^k)$ so that $a^k=\ExpSub{\xi\sim\cD}{a_{\xi}^k}$
		\STATE $h^k \eqdef a^k_{\xi^k}-a^k$
		\STATE $x^{k+1} \eqdef \left(I+\gamma A_{\xi^k}\right)^{-1}(x^k+\gamma h^k)$
		\ENDFOR
	\end{algorithmic}
	\caption{Stochastic Proximal Point Method with Operator Correction (\algname{SPPM-OC})}
	\label{alg:sppm_oc}    
\end{algorithm}

\begin{theorem}
	\label{th:sppm_gc_convergence}
	 In  Problem \eqref{eq:stoch_mon_inc_problem}, 
	let Assumptions~\ref{as:strong_mon} and \ref{as:similarity_operator} hold. 
	 Then in \algname{SPPM-OC} with any stepsize $\gamma>0$ and initial estimate $x^0\in\cX$, we have,   for every $k \geq 0$,
	\begin{equation}
		\label{eq:sppm_gc_convergence}
		\Exp{\sqn{x^{k}-x^{\star}}} \leq \left(\frac{1+\gamma^2\delta^2}{(1+\gamma\mu)^2}\right)^{k}\sqn{x^0-x^{\star}}.
	\end{equation}
	Moreover, $x^k$ converges to  $x^{\star}$, almost surely.
\end{theorem}

The contraction factor in \eqref{eq:sppm_gc_convergence} can always be made less than $1$ with $\gamma$ small enough. 
It is minimized when  $\gamma=\frac{\mu}{\delta^2}$, for which
\begin{equation*}
	\frac{1+\gamma^2\delta^2}{(1+\gamma\mu)^2}=\frac{\delta^2}{\delta^2+\mu^2}<1.
\end{equation*}
With this value of $\gamma$, the iteration complexity of \algname{SPPM-OC} to achieve $\epsilon$-accuracy for any $\epsilon>0$ is
 \begin{equation*}
\mathcal{O}\left(\left(\frac{\delta^2}{\mu^2}+1\right)\log\left(\frac{\sqn{x^{0}-x^{\star}}}{\epsilon}\right)\right).
 \end{equation*}

Thus, \algname{SPPM-OC} converges linearly to the solution $x^\star$. But it requires to select an element $a^k$ in $A(x^k)$ at every iteration, which can be costly or even impractical. Therefore, in the next section, we study another algorithm, in which this selection is performed with a small probability only.

\section{The Loopless Stochastic Variance-Reduced Proximal Point  Method (L-SVRP)}

In the optimization setting with convex \emph{differentiable} functions, the Stochastic Variance-Reduced Proximal Point Method  (SVRP) was proposed in \citet{kha23}. It was discovered  independently in \citet{tra23}, with an analysis based on the decay of the function values, which is not applicable to our setting.  This algorithm is a proximal analog of the Stochastic Variance-Reduced Gradient Method (SVRG)  \citep{joh132,zha13}, hence its name.  
 More precisely, it is a proximal analog of loopless versions of SVRG called L-SVRG \citep{hof15, kov202}. That is why we call the algorithm the Loopless Stochastic Variance-Reduced Proximal Point  Method (\algname{L-SVRP}), to emphasize its loopless nature. We introduce and study \algname{L-SVRP} in the much more general setting of set-valued monotone inclusions.

\begin{algorithm}[H]
	\begin{algorithmic}[1]
		\STATE  \textbf{Parameters:}  stepsize $\gamma>0$, initial estimates $x^0,w^0 \in\cX$, probability $p\in (0,1]$, $a^0\in A(x^0)$.
		\FOR {$k=0,1, \ldots$}
		\STATE Sample $\xi^k\sim \cD$
		\STATE Choose $a^k_{\xi^k}\in A_{\xi^k}(w^k)$ so that $\ExpSub{\xi\sim\cD}{a_{\xi}^k}=a^k$
		\STATE $h^k \eqdef a^k_{\xi^k}-a^k$
		\STATE $x^{k+1} \eqdef \left(I+\gamma A_{\xi^k}\right)^{-1}(x^k+\gamma h^k)$
		\STATE Flip a coin $\theta^k\in\{0,1\}$ with $\mathrm{Prob}(\theta^ k=1)=p$.
		\STATE $w^{k+1} \eqdef \begin{cases}
			x^{k+1}& \text{if}\ \theta^k=1\\
			w^k& \text{if}\ \theta^k=0
		\end{cases}$
		\STATE $a^{k+1} \eqdef \begin{cases}
			\text{any element in}\ A(x^{k+1})& \text{if}\ \theta^k=1\\
			a^k& \text{if}\ \theta^k=0
		\end{cases}$
		\ENDFOR
	\end{algorithmic}
	\caption{Loopless Stochastic Variance-Reduced Proximal Point Method (\algname{L-SVRP})}
	\label{alg:l_svrp}    
\end{algorithm}

\begin{theorem}
	\label{th:l_svrp_convergence}
	 In  Problem \eqref{eq:stoch_mon_inc_problem}, 
	let Assumptions~\ref{as:strong_mon} and \ref{as:similarity_operator} hold.
	 Then in \algname{L-SVRP} with any stepsize $\gamma>0$, probability  $p \in (0,1]$, and initial estimates $x^0,w^0\in\cX$, we have,   for every $k \geq 0$,
	\begin{equation}
		\label{eq:l_svrp_convergence}
		\Exp{V^{k}} \leq \max\left\{\frac{1}{1+\gamma\mu}, 1-p +\frac{\gamma\delta^2 p}{\mu\left(1+\gamma\mu\right)}\right\}^k V^0,
	\end{equation}
	where the Lyapunov function is %defined as 
	\begin{equation}
		V^{k} \eqdef \sqn{x^{k}-x^{\star}}  +\frac{\gamma\mu}{p}\sqn{w^{k} - x^\star}. 
	\end{equation}
	Moreover, $x^k$ and $w^k$ converge to  $x^{\star}$, almost surely.
\end{theorem}
The contraction factor in \eqref{eq:l_svrp_convergence} can always be made less than $1$ with $\gamma$ small enough. 
It is minimized when $\frac{1}{1+\gamma\mu}=1-p +\frac{\gamma\delta^2 p}{\mu\left(1+\gamma\mu\right)}$. This is the case for 
 \begin{equation}
 \gamma = \frac{\mu}{\delta^2+\frac{1-p}{p}\mu^2},\label{eqoptg}
 \end{equation}
for which 
 \begin{equation*}
\frac{1}{1+\gamma\mu}=1-p +\frac{\gamma\delta^2 p}{\mu\left(1+\gamma\mu\right)}=\frac{p\delta^2+(1-p)\mu^2}{p\delta^2+\mu^2}<1.
 \end{equation*}
With this value of $\gamma$, the iteration complexity of \algname{L-SVRP} to achieve $\epsilon$-accuracy for any $\epsilon>0$ is
 \begin{equation*}
\mathcal{O}\left(\left(\frac{\delta^2}{\mu^2}+\frac{1}{p}\right)\log\left(\frac{V^0}{\epsilon}\right)\right).
 \end{equation*}
The best value of $p$ depends on how much more costly it is to pick an element $a^{k}\in A(x^{k})$ than to apply the resolvent of an $A_\xi$. In any case, there is no interest in choosing $p$ larger than $\frac{\mu^2}{\delta^2}$, which is typically very small. Hence, \algname{L-SVRP} can be orders of magnitude faster than \algname{SPPM-OC}, which corresponds to the particular case of \algname{L-SVRP} with $p=1$.

In the case of minimizing a sum of $n$ differentiable functions $f_i$, i.e.\ Problem~\eqref{eq:finite_sum_mon_inc_problem} with $A_i=\nabla f_i$, with $p=\frac{1}{n}$, we recover the same  iteration complexity as in \citet{kha23}.

\section{Point-SAGA for Monotone Inclusion Problem}

\algname{Point-SAGA} is an algorithm proposed by \citet{def16} for the minimization of a sum of convex functions, using at every iteration the proximity operator of one randomly chosen function. It was also studied as a randomized primal--dual algorithm in \citet{con22rp}. 
The algorithm converges linearly when all functions are smooth and strongly convex. We introduce and study \algname{Point-SAGA} in the general setting of set-valued monotone inclusions. \algname{Point-SAGA} is an alternative to  the snapshot algorithm \algname{L-SVRP} that never requires invoking the average operator $A$. As a counterpart, 
\algname{Point-SAGA} is limited to the finite-sum problem \eqref{eq:finite_sum_mon_inc_problem}, since $n$ elements of $\cX$ are stored in a memory table.

\begin{algorithm}[H]
	\begin{algorithmic}[1]
		\STATE  \textbf{Parameters:}  stepsize $\gamma>0$, initial estimates $x^0$, $(w_i^0)_{i=1}^n \in\cX^n$, initial elements $a_i^0\in A_i(w_i^0)$ for every $i\in[n]$, $a^0\eqdef \frac{1}{n}\sum_{i=1}^n a_i^0$
		\FOR {$k=0,1, \ldots$}
		\STATE Sample $i^k \in [n]$ uniformly at random
		\STATE $h^k \eqdef  
		a_{i^k}^k - a^k$
		\STATE $x^{k+1} \eqdef \left(I + \gamma A_{i^k}\right)^{-1}\left(x^k + \gamma h^k\right)$
		\STATE ${w_j^{k+1}} \eqdef \begin{cases}
			x^{k+1}& \quad\text{for}\  \ j = i^k\\
			{ w^k_j}& \quad\text{for every}\  \ j \in [n] \backslash \{i^k\}
		\end{cases}$\ \  \ \ // not stored, defined only for the analysis
		\STATE $a_j^{k+1} \eqdef \begin{cases}
			\text{any element in}\ A_{i^k}(x^{k+1})&\ \ \text{for}\  \ j = i^k \ \ \ \ \text{// e.g. }a_{i^k}^{k+1}\eqdef\frac{1}{\gamma}(x^k-x^{k+1})+h^k\\
			a_j^k& \ \ \text{for every}\  \ j \in [n] \backslash \{i^k\}
		\end{cases}$
		\STATE $a^{k+1}\eqdef a^k + \frac{1}{n}(a_{i^k}^{k+1}-a_{i^k}^{k})$\ \ \ \ // $=\frac{1}{n}\sum_{j=1}^n a_j^{k+1}$
		\ENDFOR
	\end{algorithmic}
	\caption{\algname{Point-SAGA}}
	\label{alg:Point_SAGA}   
\end{algorithm}

\begin{theorem}
	\label{th:point_saga_convergence}
	 In  Problem \eqref{eq:finite_sum_mon_inc_problem}, 
	let Assumptions~\ref{as:strong_mon} and \ref{as:another_similarity_operator} hold.
	 Then in \algname{Point-SAGA} with any stepsize $\gamma>0$, initial estimates $x^0$, $(w_i^0)_{i=1}^n\in\cX^n$ and elements 
	$a_i^0\in A_i(w_i^0)$, we have,   for every $k \geq 0$,
	\begin{equation}
		\label{eq:point_saga_convergence}
		\Exp{V^{k}} \leq \max\left\{\frac{1}{1+\gamma\mu}, 1-\frac{1}{n} +\frac{\gamma\tilde{\delta}^2}{n\mu(1+\gamma\mu)}\right\}^k V^0,
	\end{equation}
	where the Lyapunov function is 
	\begin{equation}
		V^{k} \eqdef \sqn{x^k-x^{\star}}  +\gamma\mu\sum_{i=1}^n\sqn{w^k_i - x^{\star}}. 
	\end{equation}
	Moreover, $x^k$ and all $w^k_i$ converge to  $x^{\star}$, almost surely.
\end{theorem}
The contraction factor in \eqref{eq:point_saga_convergence} can always be made less than $1$ with $\gamma$ small enough. 
It is minimized when $\frac{1}{1+\gamma\mu}=1-\frac{1}{n} +\frac{\gamma\tilde{\delta}^2 }{n\mu(1+\gamma\mu)}$. This is the case for 
 \begin{equation*}
 \gamma = \frac{\mu}{\tilde{\delta}^2+(n-1)\mu^2},
 \end{equation*}
for which 
 \begin{equation*}
\frac{1}{1+\gamma\mu}=1-\frac{1}{n} +\frac{\gamma\tilde{\delta}^2 }{n\mu(1+\gamma\mu)}=\frac{\tilde{\delta}^2+(n-1)\mu^2}{\tilde{\delta}^2+n\mu^2}<1.
 \end{equation*}
With this value of $\gamma$, the iteration complexity of \algname{Point-SAGA} to achieve $\epsilon$-accuracy for any $\epsilon>0$ is
 \begin{equation*}
\mathcal{O}\left(\left(\frac{\delta^2}{\mu^2}+n\right)\log\left(\frac{V^0}{\epsilon}\right)\right).
 \end{equation*}
 
To the best of our knowledge, the analysis of \algname{Point-SAGA} under a similarity assumption is new, even in the particular case of minimizing convex  functions.

\section{Experiments}

\begin{figure}[t]
	\centering
	\includegraphics[width=0.32\textwidth]{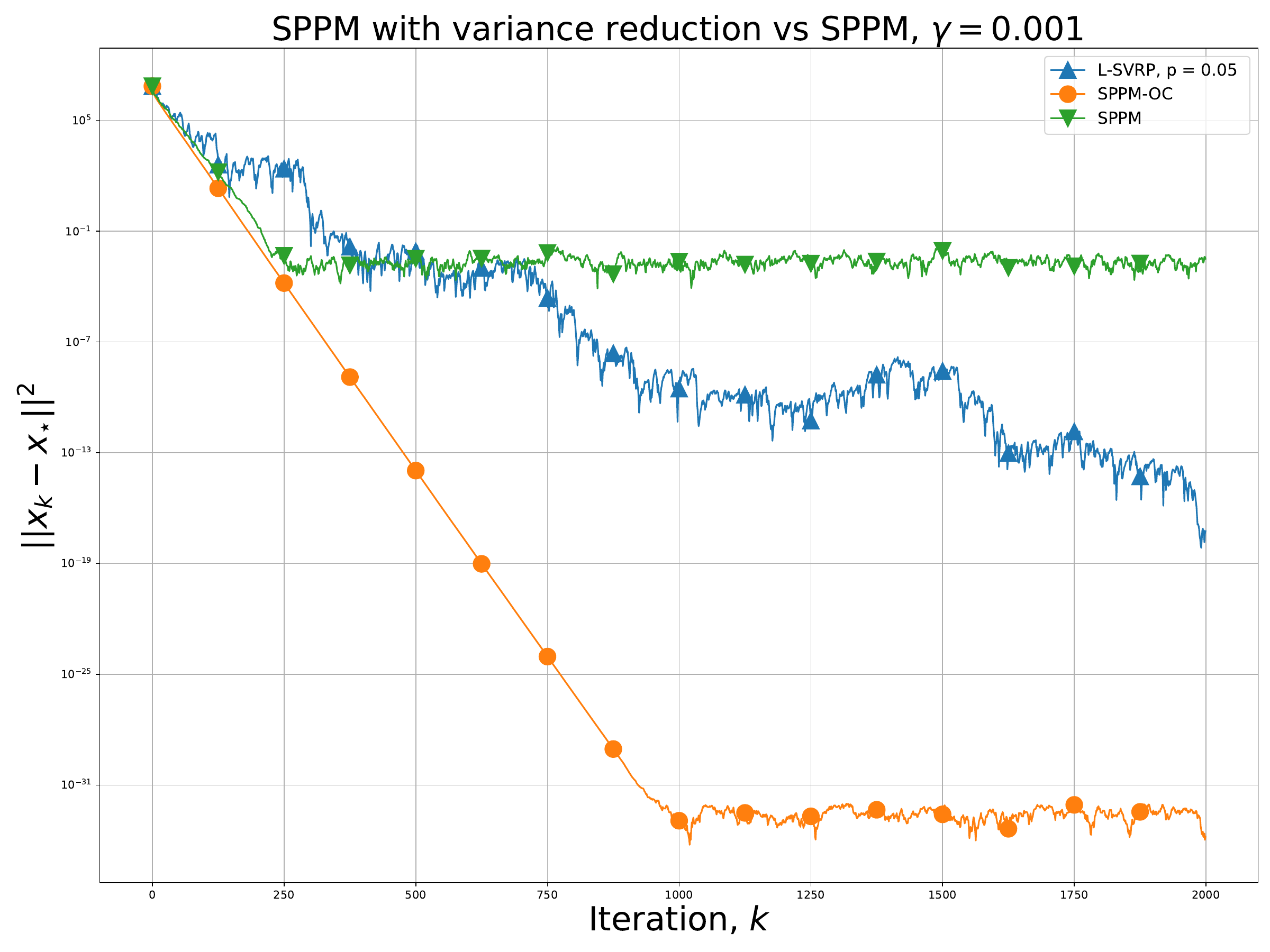}
	\includegraphics[width=0.32\textwidth]{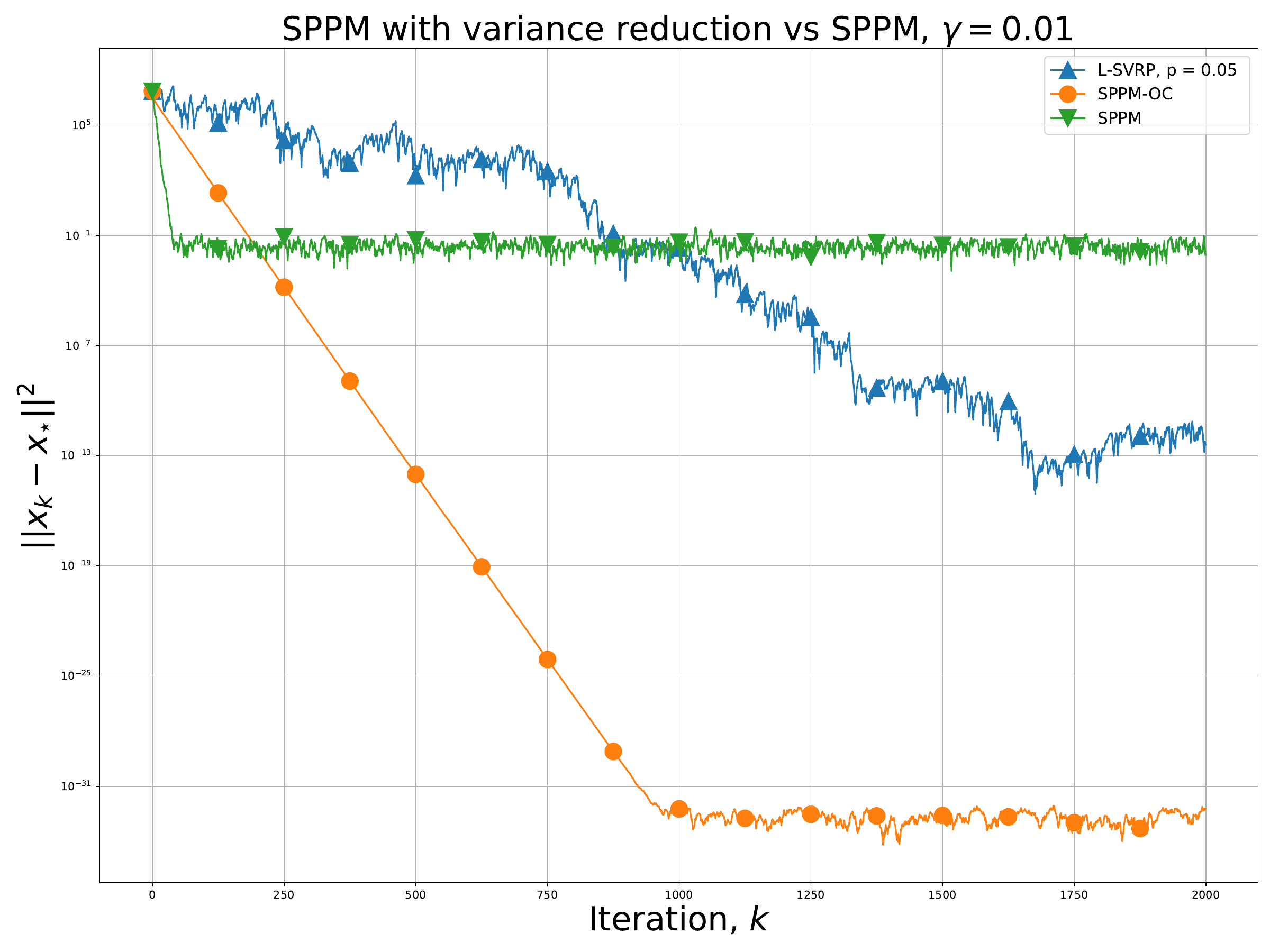}
	\includegraphics[width=0.32\textwidth]{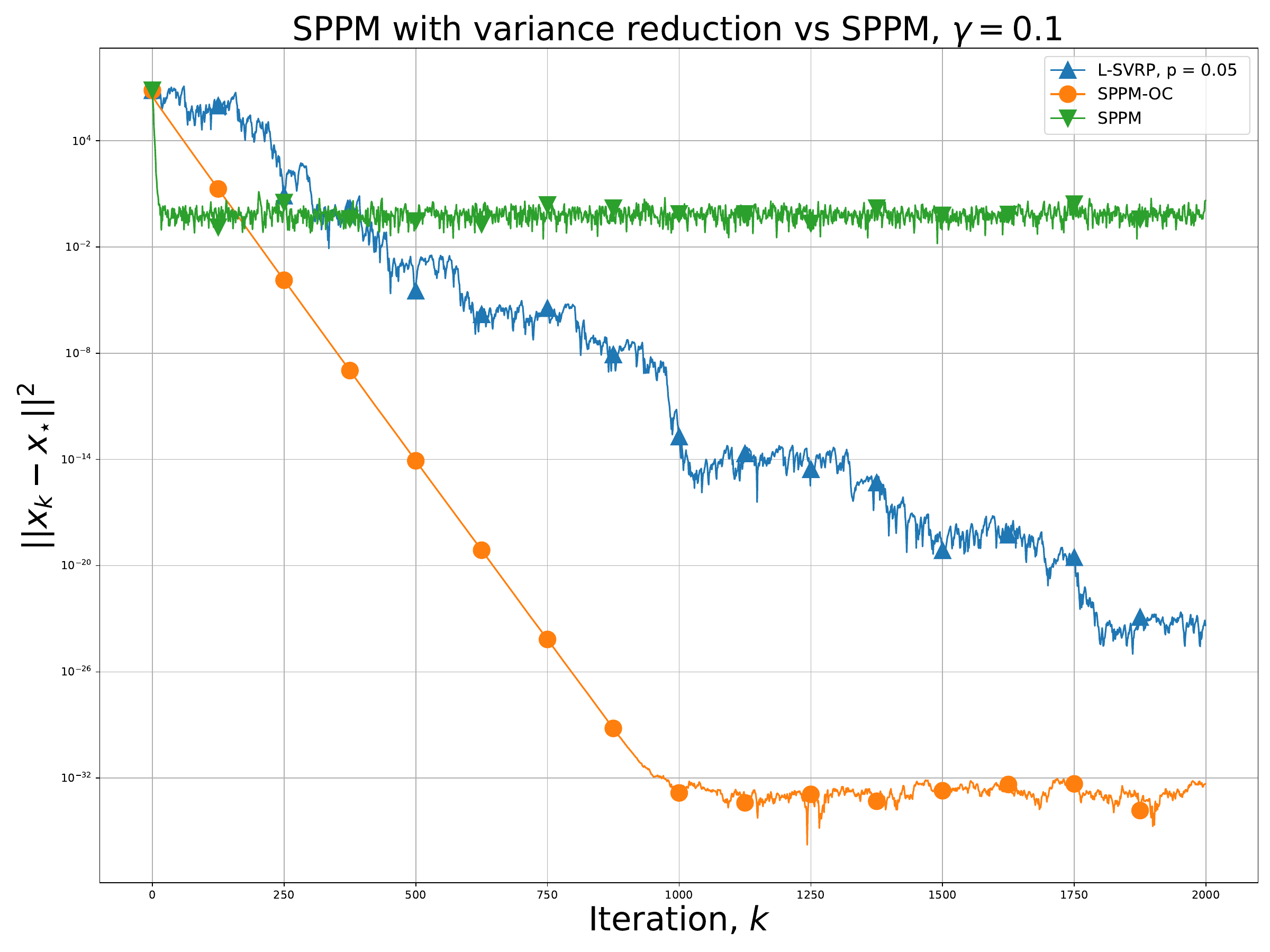}
	\caption{Performance comparison of \algname{SPPM} with, from left to right, $\gamma =10^{-3}, 10^{-2}, 10^{-1}$,  \algname{SPPM-OC}, \algname{L-SVRP}. \algname{SPPM-OC} and \algname{L-SVRP} have the same parameter values in the 3 plots, the differences are only due to randomness.}
	\label{fig:exp_1}
\end{figure}

We perform numerical experiments for the saddle-point problem 
\begin{equation*}
	\min_{y\in \R^{d_y}}\max_{z \in \R^{d_z} } \frac{1}{n}\sum^n_{i}f_i(y,z),
\end{equation*}
for some vector dimensions $d_y\geq 1$ and $d_z\geq 1$, 
where each $f_i$ is a strongly convex--strongly concave function defined as
\begin{equation*}
	f_i:(y,z) \mapsto \frac{1}{2}\abr{y, M_{i}y} +\abr{b_{i},y} + \abr{z,Q_i y} -\abr{c_{i}, z} - \frac{1}{2}\abr{z, N_{i} z},
\end{equation*} 
with the following parameters: 
\begin{itemize}
	\item Each matrix $M_i \in \R^{d_y\times d_y}$ and $N_i  \in \R^{d_z\times d_z}$ is generated randomly with apriori selected eigenvalues $\lambda_l(M_i)= 10^l$ and  $\lambda_j(N_i)= 10^j$ respectively, where $l \in \{0, 1, \dots, d_y -1\}$ and $j \in \{0, 1, \dots, d_z -1\}$;
	\item  The vectors $b_i \in \R^{d_y} $ and $c_i \in \R^{d_z}$ are sampled from normal distributions $\cN(1, 5\cdot I_{d_y})$ and $\cN(1, 5 \cdot I_{d_z})$ respectively;
	\item Every element of the matrix $Q_i$ is sampled from the standard normal distribution, then each column is normalized  to have a full-rank matrix. 
\end{itemize}
To formulate the problem as Problem~\eqref{eq:finite_sum_mon_inc_problem}, we define $x\eqdef (y,z)$ and the single-valued monotone operators
$A_i :x\mapsto  (\nabla_y^{\top} f_i(y,z), -\nabla_z^{\top} f_i(y,z))^{\top} $; that is,
\begin{equation*}
	A_i(x) =\begin{pNiceArray}{c|c}
		M_i & Q_i^{\top} \\
		\hline
		-Q_i &  N_i
	\end{pNiceArray} x + \begin{pmatrix} b_i \\ c_i \end{pmatrix} = \mathbb{B}_ix + r_i.
\end{equation*}
We take $n = 200$, $d_y =3$, $d_z = 4$. Each operator $A_i$ is $1$-strongly monotone and $L$-Lipschitz-continuous with $L = 1000$.

We compute the similarity constant $\delta$ as follows. By Assumption~\ref{as:similarity_operator}, we have 
\begin{equation*}
	 \frac{1}{n}\sum_{i=1}^n\sqn{A_i(x) - A(x) - A_i(x^{\star}) + A(x^{\star})} \leq \delta^2  \sqn{x-x^{\star}},
\end{equation*}
Plugging in the expression for $A_i (x)= \mathbb{B}_i x + r_i$, we obtain 
\begin{align*}
	 \frac{1}{n}\sum_{i=1}^n\sqn{A_i(x) - A(x) - A_i(x^{\star}) + A(x^{\star})}   &= \frac{1}{n}\sum_{i=1}^n\sqn{\mathbb{B}_i x - \frac{1}{n}\sum_{j=1}^n \mathbb{B}_j x - \mathbb{B}_i x^{\star} + \frac{1}{n}\sum_{j=1}^n \mathbb{B}_j  x^{\star} } \\
	 &\leq \frac{1}{n}\sum_{i=1}^n\sqn{\mathbb{B}_i  - \frac{1}{n}\sum_{j=1}^n \mathbb{B}_j } \sqn{x-x^{\star}}.
\end{align*}
Thus we have the simple and easy to compute upper bound
\begin{equation}
	\delta \leq  \frac{1}{n}\sum_{i=1}^n\sqn{\mathbb{B}_i  - \frac{1}{n}\sum_{j=1}^n \mathbb{B}_j } \approx 26.5.
\end{equation}
As we can see, $\delta<\!\!<L$.

\begin{figure}[t]
	\centering
	\includegraphics[width=0.49\textwidth]{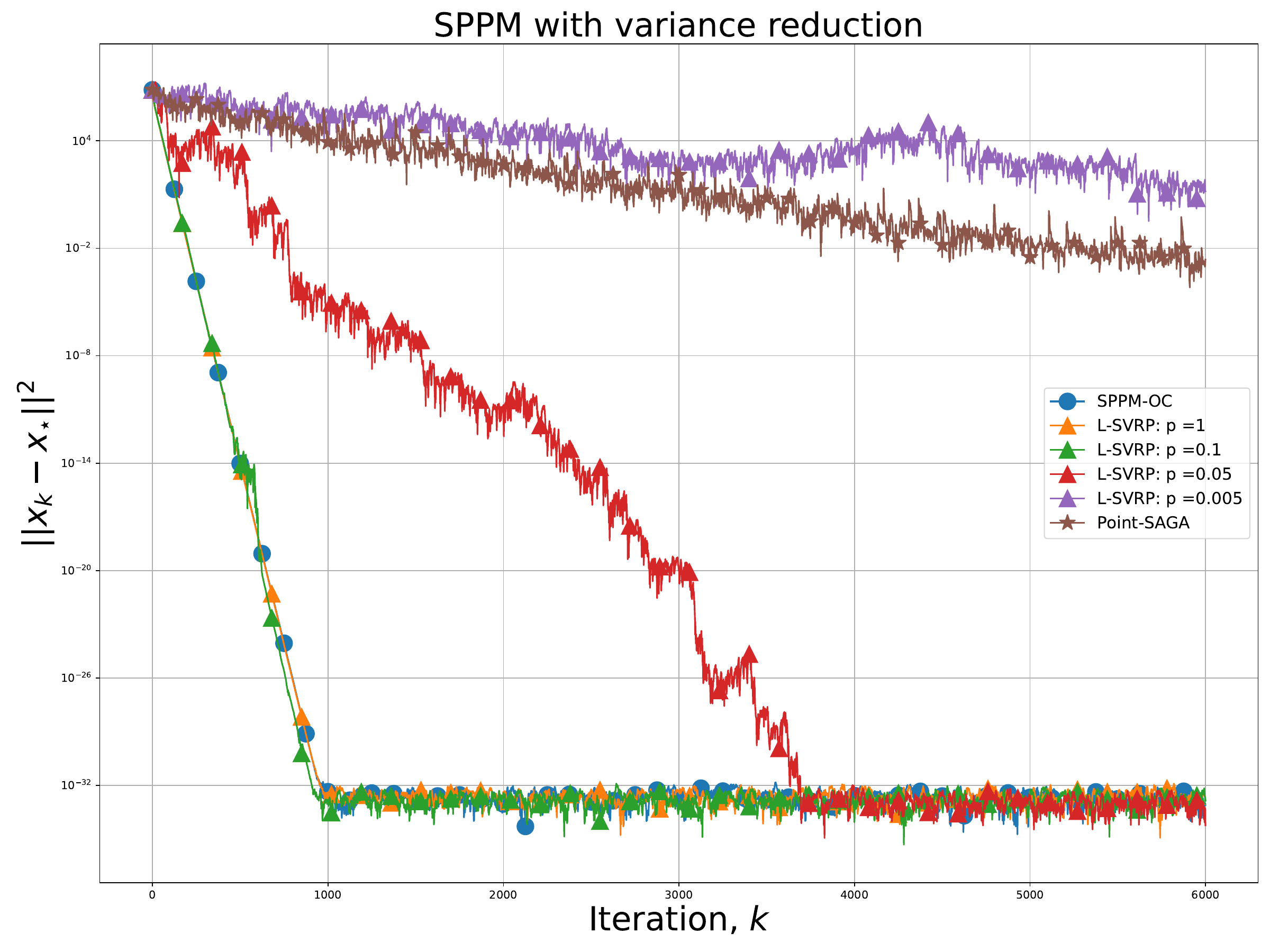}
	\includegraphics[width=0.49\textwidth]{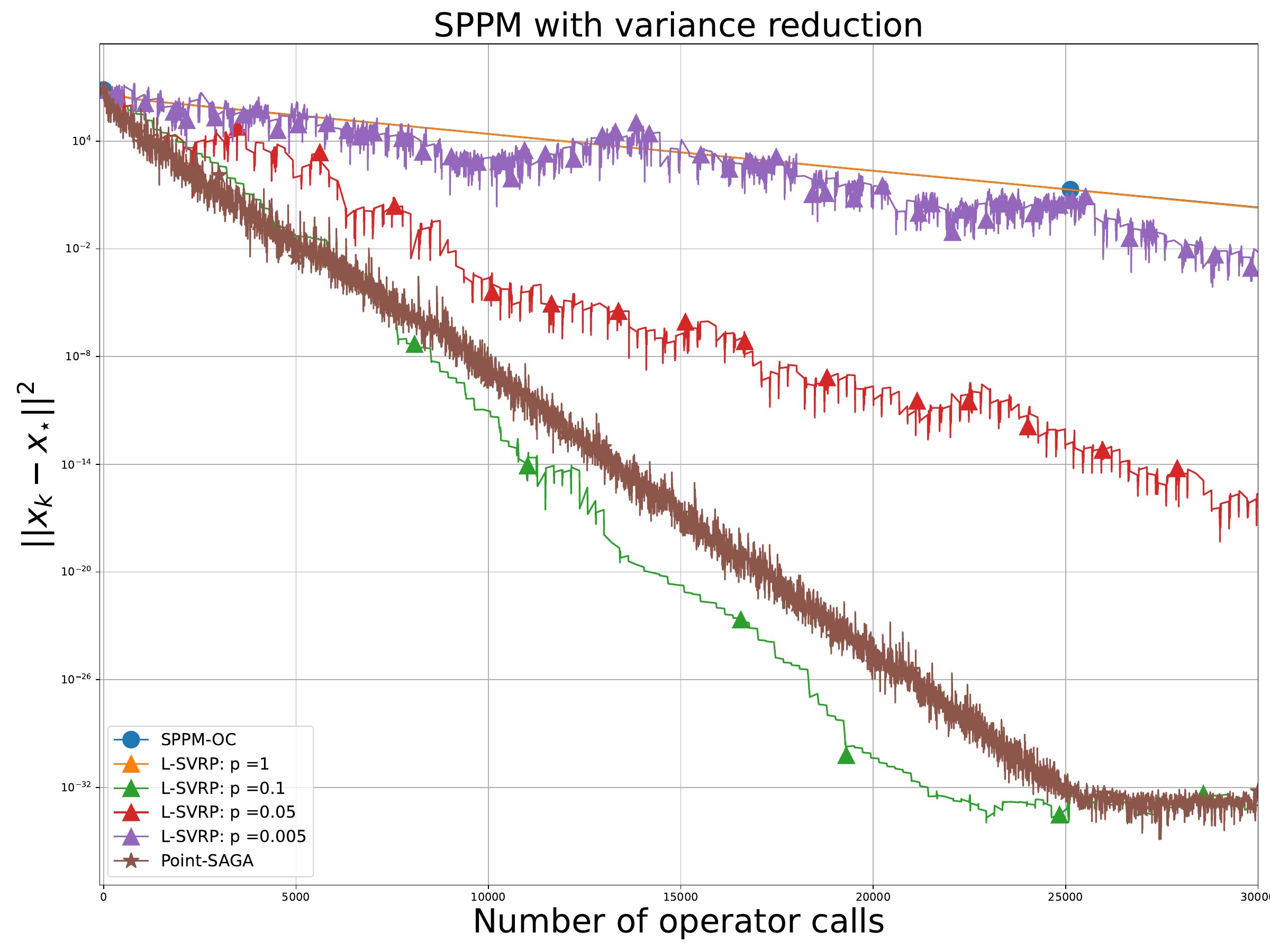}
	\caption{Performance comparison of  \algname{SPPM-OC},  \algname{L-SVRP} with different values  of  $p =1, 0.1, 0.05, \nicefrac{1}{n} =0.005$, and  \algname{Point-SAGA}. The theoretically optimal value of $\gamma$ is chosen in all cases. The error is shown with respect to the number of iterations on the left, and  the  number of operator calls on the right.
	}
	\label{fig:exp_2}
\end{figure}

In Figure~\ref{fig:exp_1}, we compare  \algname{SPPM} with 3 different values of $\gamma$, \algname{SPPM-OC} with the theoretically optimal value $\gamma  = \frac{\mu}{\delta^2} \approx 10^{-3}$,  \algname{L-SVRP} with $p=0.05$ and the theoretically optimal value of $\gamma$ in \eqref{eqoptg}. As predicted by the theory,  \algname{SPPM} converges only to a neighborhood of the solution, whose size is larger if $\gamma$ is larger.  \algname{SPPM-OC} is faster than  \algname{L-SVRP}, but its per-iteration cost is much higher, as we detail in Figure~\ref{fig:exp_2}.

In Figure~\ref{fig:exp_2}, we compare  the variance-reduced algorithms \algname{SPPM-OC},   \algname{L-SVRP} with different values of $p$,   and \algname{Point-SAGA}. 
 \algname{SPPM-OC} and \algname{L-SVRP} with $p=1$ are identical. 
We show convergence with respect to the number of operator calls, counting 1 for  a call to an $A_\xi$ or its resolvent, and $n$ for a call to $A$ in  \algname{L-SVRP}. As a result,  \algname{L-SVRP} with $p=0.1$ and  \algname{Point-SAGA} perform best. We should keep in mind that  \algname{L-SVRP} does a full pass over the $n$ operators with a small probability, whereas  \algname{Point-SAGA} requires memory storage of size $n$ times the dimension of $x$. Thus, the best algorithm depends on the problem at hand.

\section*{Acknowledgement}
This work was supported by the SDAIA-KAUST Center of Excellence in Data Science and Artificial Intelligence (SDAIA-KAUST AI).

\bibliographystyle{icml2024}
\bibliography{IEEEabrv,biblio2}

\newpage
\appendix

\noindent{\huge\textbf{Appendix}}

\section{Proof of Lemma~\ref{fact:non-expansiveness}}
Let $(x,y)\in\cX^2$. From the definition of the resolvent, we have $x-x^+ \in  B(x^+)$ and $y-y^+ \in B(y^+)$. Then it follows from  \eqref{eq:strong_mon} that
\begin{equation*}
\mu\sqn{x^+-y^+}\leq \abr{(x-x^+)-(y-y^+), x^+-y^+}=\abr{x-y, x^+-y^+}-\sqn{x^+-y^+}.
\end{equation*}
Therefore
\begin{equation*}
(1+\mu)\sqn{x^+-y^+}\leq \abr{x-y, x^+-y^+} \leq \|x-y\|\| x^+-y^+\|,
\end{equation*}
so that 
\begin{equation*}
(1+\mu)\|x^+-y^+\| \leq \|x-y\|.
\end{equation*}

\section{Proof of Theorem~\ref{th:sppm_convergence}}
Let $k\geq 0$.
    We have 
    \begin{equation}
    \label{eq16}
    x^{\star} = \left(I+\gamma A_{\xi^k}\right)^{-1}\left(x^{\star} + \gamma a_{\xi^k}^{\star} \right),
    \end{equation}
    so that
    \begin{eqnarray*}
        \sqn{x^{k+1}-x^{\star}} &=& \sqn{\left(I+\gamma A_{\xi^k}\right)^{-1}(x^k)-\left(I+\gamma A_{\xi^k}\right)^{-1}\left(x^{\star} + \gamma a_{\xi^k}^\star \right)}\\
        &\overset{\eqref{eq:non-expansiveness}}\leq&\frac{1}{(1+\gamma\mu)^2}\sqn{x^k - x^{\star} - \gamma a_{\xi^k}^\star}\\
        &=& \frac{1}{(1+\gamma\mu)^2}\left(\sqn{x^k - x^{\star}} -2\gamma\abr{a_{\xi^k}^\star, x^k - x^{\star}} + \gamma^2\sqn{a_{\xi^k}^\star}\right).
    \end{eqnarray*}
    We denote by $\cF^k$ the $\sigma$-algebra generated by the collection of random variables $(x^0,\ldots, x^k)$.
Taking the expectation conditionally on $\cF^k$, we have, using the fact that $\ExpSub{\xi\sim\cD}{a_{\xi}^\star}=0$,
\begin{eqnarray*}
    \ExpCond{\sqn{x^{k+1}-x^{\star}}}{\cF^k} &\leq& \frac{1}{(1+\gamma\mu)^2}\left(\sqn{x^k - x^{\star}} -2\gamma\abr{\underbrace{\ExpCond{a_{\xi^k}^\star}{\cF^k}}_0, x^k - x^{\star}}\right)\\
    &&\quad{}+ \frac{\gamma^2}{(1+\gamma\mu)^2}\ExpCond{\sqn{a_{\xi^k}^\star}}{\cF^k}\\
    &\overset{\eqref{eq:sigma_star}}{=}& \frac{1}{(1+\gamma\mu)^2}\sqn{x^k - x^{\star}} + \frac{\gamma^2\sigma^2_{\star}}{(1+\gamma\mu)^2}.
\end{eqnarray*} 
By unrolling the recursion, we obtain
\begin{eqnarray*}
    \Exp{\sqn{x^{k}-x^{\star}}} &\leq& \left(\frac{1}{1+\gamma\mu}\right)^{2k} \sqn{x^0 - x^{\star}} +\sum^{k-1}_{l=0} \left(\frac{1}{1+\gamma\mu}\right)^{2l} \frac{\gamma^2\sigma^2_{\star}}{(1+\gamma\mu)^2}\\
     &=& \left(\frac{1}{1+\gamma\mu}\right)^{2k} \sqn{x^0 - x^{\star}} + \frac{(1+\gamma\mu)^2-(1+\gamma\mu)^{2(1-k)}}{(1+\gamma\mu)^2-1}\frac{\gamma^2\sigma^2_{\star}}{(1+\gamma\mu)^2}\\
          &=& \left(\frac{1}{1+\gamma\mu}\right)^{2k} \sqn{x^0 - x^{\star}} + \frac{1-(1+\gamma\mu)^{-2k}}{(1+\gamma\mu)^2-1}\gamma^2\sigma^2_{\star}\\
                    &\leq& \left(\frac{1}{1+\gamma\mu}\right)^{2k} \sqn{x^0 - x^{\star}} + \frac{1}{(1+\gamma\mu)^2-1}\gamma^2\sigma^2_{\star}\\
             &=& \left(\frac{1}{1+\gamma\mu}\right)^{2k} \sqn{x^0 - x^{\star}} + \frac{\gamma\sigma^2_{\star}}{2\mu+\gamma\mu^2}.
\end{eqnarray*}

\section{Proof of Theorem~\ref{th:sppm_gc_convergence}}
%\begin{proof}
    For every $\xi\sim\cD$, let $a_\xi^\star \in A_{\xi}(x^{\star})$, such  that $\ExpSub{\xi\sim\cD}{a_\xi^\star}=0$ and Assumption \ref{as:similarity_operator} holds at $x^k$ with these elements.
    Let $k\geq 0$. Using \eqref{eq16}, 
 	we have 
	\begin{eqnarray*}
		\sqn{x^{k+1}-x^{\star}} &=& \sqn{\left(I+\gamma A_{\xi^k}\right)^{-1}(x^k +\gamma h^k)-\left(I+\gamma A_{\xi^k}\right)^{-1}\left(x^{\star} + \gamma a_{\xi^k}^\star \right)}\\
		&\overset{\text{Lemma~\ref{fact:non-expansiveness}}}\leq&\frac{1}{(1+\gamma\mu)^2}\sqn{x^k - x^{\star} + \gamma h^k - \gamma a_{\xi^k}^\star }\\
		&=& \frac{1}{(1+\gamma\mu)^2}\left(\sqn{x^k - x^{\star}} +2\gamma\abr{h^k - a_{\xi^k}^\star, x^k - x^{\star}} + \gamma^2\sqn{h^k -a_{\xi^k}^\star }\right).
	\end{eqnarray*}
	 We denote by $\cF^k$ the $\sigma$-algebra generated by the collection of random variables $(x^l,a^l, a_{\xi^l}^l)_{l=0}^k$.
Taking the expectation conditionally on $\cF^k$, we have, using the fact that $\ExpSub{\xi\sim\cD}{a_{\xi}^\star}=0$,
	\begin{eqnarray}
		\ExpCond{\sqn{x^{k+1}-x^{\star}}}{\cF^k} &\leq& \frac{1}{(1+\gamma\mu)^2}\sqn{x^k - x^{\star}} \notag\\
		&&+\frac{2\gamma}{(1+\gamma\mu)^2}\abr{\underbrace{\ExpCond{a^k_{\xi^k}-a^k -a_{\xi^k}^\star}{\cF^k}}_0, x^k - x^{\star}}\notag\\
		&&+ \frac{\gamma^2}{(1+\gamma\mu)^2}\ExpCond{\sqn{a^k_{\xi^k}-a^k - a_{\xi^k}^\star}}{\cF^k}\notag\\
		&\overset{\eqref{eq:similarity_operator}}{\leq}& \frac{1}{(1+\gamma\mu)^2}\sqn{x^k - x^{\star}} + \frac{\gamma^2\delta^2}{(1+\gamma\mu)^2}\sqn{x^k - x^{\star}}\notag\\
		&=& \frac{1+\gamma^2\delta^2}{(1+\gamma\mu)^2}\sqn{x^k - x^{\star}}.\label{eqrec2b}
	\end{eqnarray} 
	By unrolling the recursion, we obtain the desired result.
%\end{proof}
Moreover, using classical results on supermartingale convergence \citep[Proposition A.4.5]{ber15}, it follows from \eqref{eqrec2b} that 
$\sqn{x^k - x^{\star}}\rightarrow 0$ almost surely.

\section{Proof of Theorem~\ref{th:l_svrp_convergence}}
    For every $\xi\sim\cD$, let $a_\xi^\star \in A_{\xi}(x^{\star})$, such  that $\ExpSub{\xi\sim\cD}{a_\xi^\star}=0$ and Assumption \ref{as:similarity_operator} holds at $x^k$ with these elements. 
	Let $k\geq 0$. Using \eqref{eq16}, 
	we have 
	\begin{eqnarray*}
		\sqn{x^{k+1}-x^{\star}} &=& \sqn{\left(I+\gamma A_{\xi^k}\right)^{-1}(x^k +\gamma h^k)-\left(I+\gamma A_{\xi^k}\right)^{-1}\left(x^{\star} + \gamma a_{\xi^k}^\star \right)}\\
		&\overset{\text{Lemma~\ref{fact:non-expansiveness}}}\leq&\frac{1}{(1+\gamma\mu)^2}\sqn{x^k - x^{\star} + \gamma h^k - \gamma a_{\xi^k}^\star }\\
		&=& \frac{1}{(1+\gamma\mu)^2}\left(\sqn{x^k - x^{\star}} +2\gamma\abr{h^k - a_{\xi^k}^\star, x^k - x^{\star}} + \gamma^2\sqn{h^k -a_{\xi^k}^\star }\right).
	\end{eqnarray*}
	 We denote by $\cF^k$ the $\sigma$-algebra generated by the collection of random variables $(x^l,w^l,a^l, a_{\xi^l}^l)_{l=0}^k$.
Taking the expectation conditionally on $\cF^k$, we have, using the fact that $\ExpSub{\xi\sim\cD}{a_{\xi}^\star}=0$,
	\begin{eqnarray}
		\ExpCond{\sqn{x^{k+1}-x^{\star}}}{\cF^k} &\leq& \frac{1}{(1+\gamma\mu)^2}\sqn{x^k - x^{\star}}\notag\\
		&&+\frac{2\gamma}{(1+\gamma\mu)^2}\abr{\underbrace{\ExpCond{a^k_{\xi^k}-a^k -a_{\xi^k}^\star}{\cF^k}}_0, x^k - x^{\star}}\notag\\
		&&+ \frac{\gamma^2}{(1+\gamma\mu)^2}\ExpCond{\sqn{a^k_{\xi^k}-a^k - a_{\xi^k}^\star}}{\cF^k}\notag\\
		&\overset{\eqref{eq:similarity_operator}}{\leq}& \frac{1}{(1+\gamma\mu)^2}\sqn{x^k - x^{\star}} + \frac{\gamma^2\delta^2}{(1+\gamma\mu)^2}\sqn{w^k - x^{\star}}.\label{eq:uyweiudncxjk}
	\end{eqnarray} 
Moreover,
\begin{eqnarray*}
	\ExpCond{\sqn{w^{k+1} - x^{\star}}}{\cF^k} &=& (1-p)\|w^k-x^{\star}\|^2 + p\ExpCond{
		\sqn{x^{k+1} -x^{\star}}}{\cF^k}.
\end{eqnarray*}
Let $\alpha \eqdef \frac{\gamma\mu}{p}$. 
Combining the two previous inequalities and using the Lyapunov function $V^{k+1} \eqdef \sqn{x^{k+1}-x^{\star}}  +\alpha\sqn{w^{k+1} - x^{\star}}$, we obtain
\begin{eqnarray}
	\ExpCond{V^{k+1}}{\cF^k} &\leq& \ExpCond{\sqn{x^{k+1} - x^{\star}}}{\cF^k} + \alpha p\ExpCond{\sqn{x^{k+1} -x^{\star}}}{\cF^k}\notag\\
	&& + (1-p)\alpha\sqn{w^k-x^{\star}}\notag\\
	&=& (1+\alpha p)\ExpCond{\sqn{x^{k+1} - x^{\star}}}{\cF^k} +  (1-p)\alpha\sqn{w^k-x^{\star}}\notag\\
	&\overset{\eqref{eq:uyweiudncxjk}}{\leq}& \frac{1+\alpha p}{(1+\gamma\mu)^2}\sqn{x^k - x^{\star}} + \frac{(1+\alpha p)\gamma^2\delta^2}{(1+\gamma\mu)^2}\sqn{w^k - x^{\star}} +  (1-p)\alpha\sqn{w^k-x^{\star}}\notag\\
	&=& \frac{1+\alpha p}{(1+\gamma\mu)^2}\sqn{x^k - x^{\star}} + \left(1-p + \frac{(1+\alpha p)\gamma^2\delta^2}{\alpha(1+\gamma\mu)^2}\right)\alpha\sqn{w^k-x^{\star}}\notag\\
	&\overset{\alpha = \frac{\gamma\mu}{p}}{\leq}& \max\left\{\frac{1}{1+\gamma\mu}, 1-p+\frac{\gamma\delta^2 p}{\mu(1+\gamma\mu)}\right\} V^k. \label{eqrec2c} 
\end{eqnarray}
By unrolling the recursion, we obtain the desired result. Moreover, using classical results on supermartingale convergence \citep[Proposition A.4.5]{ber15}, it follows from \eqref{eqrec2c} that 
$V^k\rightarrow 0$ almost surely. 

\section{Proof of Theorem~\ref{th:point_saga_convergence}}

For every $i\in[n]$, let $a_i^\star \in A_{i}(x^{\star})$, such  that $\frac{1}{n}\sum_{i=1}^n a_i^\star=0$ and Assumption \ref{as:another_similarity_operator} holds at the $(w_i^k)_{i=1}^n$ with these elements.
Let $k\geq 0$. We have
   \begin{equation*}
    x^{\star} = \left(I+\gamma A_{i^k}\right)^{-1}\left(x^{\star} + \gamma a_{i^k}^{\star} \right),
    \end{equation*}
so that 
	\begin{eqnarray*}
		\sqn{x^{k+1}-x^{\star}} &=& \sqn{\left(I+\gamma A_{i^k}\right)^{-1}(x^k +\gamma h^k)-\left(I+\gamma A_{i^k}\right)^{-1}\left(x^{\star} + \gamma a_{i^k}^\star \right)}\\
		&\overset{\text{Lemma~\ref{fact:non-expansiveness}}}\leq&\frac{1}{(1+\gamma\mu)^2}\sqn{x^k - x^{\star} + \gamma h^k - \gamma a_{i^k}^\star }\\
		&=& \frac{1}{(1+\gamma\mu)^2}\left(\sqn{x^k - x^{\star}} +2\gamma\abr{h^k - a_{i^k}^\star, x^k - x^{\star}} + \gamma^2\sqn{h^k -a_{i^k}^\star }\right).
	\end{eqnarray*}
	 We denote by $\cF^k$ the $\sigma$-algebra generated by the collection of random variables $\big(x^l,(w_i^l)_{i=1}^n,(a_i^l)_{i=1}^n\big)_{l=0}^k$.
Taking the expectation conditionally on $\cF^k$, we have
\begin{eqnarray}
	\ExpCond{\sqn{x^{k+1}-x^{\star}}}{\cF^k} &\leq& \frac{1}{(1+\gamma\mu)^2}\sqn{x^k - x^{\star}}\notag\\
	&& + \frac{\gamma^2}{(1+\gamma\mu)^2}\ExpCond{\sqn{a_{i^k}^k - a^k- a_{i^k}^{\star}}}{\cF^k}\notag\\
	&&+ \frac{2\gamma}{(1+\gamma\mu)^2}\abr{\underbrace{\ExpCond{a_{i^k}^k - a^k- a_{i^k}^{\star}}{\cF^k}}_0, x^k - x^{\star}}\notag\\
	&=& \frac{1}{(1+\gamma\mu)^2}\sqn{x^k - x^{\star}} + \frac{\gamma^2}{n(1+\gamma\mu)^2}\sum^n_{i=1}\sqn{a_{i}^k - a^k- a_{i}^{\star}}\notag\\
	&\overset{\eqref{eq:another_similarity_operator}}{\leq}& \frac{1}{(1+\gamma\mu)^2}\sqn{x^k - x^{\star}} + \frac{\gamma^2\tilde{\delta}^2}{n(1+\gamma\mu)^2}\sum^n_{i=1}\sqn{w^k_i - x^{\star}} \label{eq:wopqwpqsdjsj}.
\end{eqnarray} 
Moreover,
\begin{eqnarray*}
	\frac{1}{n}\sum^n_{i=1}\ExpCond{\sqn{w_i^{k+1} - x^{\star}}}{\cF^k} &=& \left(1-\frac{1}{n}\right)\frac{1}{n}\sum^n_{i=1}\sqn{w^k_i-x^{\star}} + \frac{1}{n}\ExpCond{\sqn{x^{k+1} -x^{\star}}}{\cF^k}.
\end{eqnarray*}
Let $\alpha \eqdef n\gamma\mu$. 
Combining the two previous inequalities and using the Lyapunov function $V^{k+1} \eqdef \sqn{x^{k+1}-x^{\star}}  +\frac{\alpha}{n}\sum^n_{i=1}\sqn{w_i^{k+1} - x^{\star}}$, we obtain
\begin{eqnarray}
	\ExpCond{V^{k+1}}{\cF^k} &\leq& \ExpCond{\sqn{x^{k+1} - x^{\star}}}{\cF^k} + \frac{\alpha}{n} \ExpCond{\sqn{x^{k+1} -x^{\star}}}{\cF^k}\notag\\
	&& + \left(1-\frac{1}{n}\right)\frac{\alpha}{n}\sum^n_{i=1}\sqn{w^k_i-x^{\star}}\notag\\
	&=& \left(1+\frac{\alpha}{n} \right)\ExpCond{\sqn{x^{k+1} - x^{\star}}}{\cF^k} +  \left(1-\frac{1}{n}\right)\frac{\alpha}{n}\sum^n_{i=1}\sqn{w^k_i-x^{\star}}\notag\\
	&\overset{\eqref{eq:wopqwpqsdjsj}}{\leq}& \frac{1+\nicefrac{\alpha}{n}}{(1+\gamma\mu)^2}\sqn{x^k - x^{\star}} + \frac{(1+\nicefrac{\alpha}{n})\gamma^2\tilde\delta^2}{(1+\gamma\mu)^2}\frac{1}{n}\sum^n_{i=1}\sqn{w^k_i - x^{\star}}\notag\\
	&& + \left(1-\frac{1}{n}\right)\frac{\alpha}{n}\sum^n_{i=1}\sqn{w^k_i-x^{\star}}\notag\\
	&=& \frac{1+\nicefrac{\alpha}{n}}{(1+\gamma\mu)^2}\sqn{x^k - x^{\star}} + \left(1-\frac{1}{n} + \frac{(1+\nicefrac{\alpha}{n})\gamma^2\tilde\delta^2}{\alpha(1+\gamma\mu)^2}\right)\frac{\alpha}{n}\sum^n_{i=1}\sqn{w^k_i-x^{\star}}\notag\\
	&\overset{\alpha = n\gamma\mu}{\leq}& \max\left\{\frac{1}{1+\gamma\mu}, 1-\frac{1}{n}+\frac{\gamma\tilde{\delta}^2 }{n\mu(1+\gamma\mu)}\right\} V^k.\label{eqrec2d}
\end{eqnarray}
By unrolling the recursion, we obtain the desired result. Moreover, using classical results on supermartingale convergence \citep[Proposition A.4.5]{ber15}, it follows from \eqref{eqrec2d} that 
$V^k\rightarrow 0$ almost surely.

\end{document}